\documentclass[a4paper,11pt]{article}
 \usepackage[T1]{fontenc}
 \usepackage[normalem]{ulem}
 \usepackage{verbatim}
 \usepackage{graphicx}
 \usepackage{amsmath}
 \usepackage{amssymb}
\usepackage{geometry}
\usepackage{amsthm}
\usepackage[utf8]{inputenc}

\newtheorem{proposition}{Proposition}[section]

\newtheorem{theorem}{Theorem}[section]

\newtheorem{e-proposition}[theorem]{Proposition}
\newtheorem{corollary}[theorem]{Corollary}
\newtheorem{e-definition}[theorem]{Definition\rm}

\def\cC{\mathcal{C}}

\def\N{\mathbb{N}}

\def\R{\mathbb{R}}

\title{Global existence of bounded weak solutions to degenerate cross-diffusion equations in moving domain}
\author{Athmane Bakhta, Virginie Ehrlacher}
\date{\today}

\begin{document}

\maketitle

\section{Introduction}
\label{sec:introduction}
This note focuses on some issues for the analysis of a system of degenerate cross-diffusion partial differential equations (PDEs). 
This family of models are encountered in a wide variety of contexts, such as population dynamics 
~\cite{JungelZamponi}, biology~\cite{Painter2,Painter1}, chemistry or materials science. 
The application we have in mind here is the modeling of the evolution of the concentration of chemical species composing a crystalline solid. The functions, that are the solutions 
of the system of PDEs of interest, represent the local densities of the different components of the material, and thus should be nonnegative, bounded and satisfy some volumic constraints which will be made precise 
later in the note. These systems are useful for instance for the prediction of the chemical composition of thin solid films grown by Chemical Vapor Deposition (CVD)~\cite{CVD}. In this process, 
a solid wafer is exposed to gaseous precursors, corresponding to the different species entering the composition of the film, 
which react or decompose on the substrate surface to produce the desired deposit. This process generally occurs at high temperature and takes several hours, 
so that the diffusion of the different atomic species within the bulk of the 
solid has to be taken into account in addition to the evolution of the surface of the film.

\medskip

More precisely, let $d\in \N^*$ (the space dimension), $T>0$ (the duration of the process), and for all $t\in[0,T]$, let $\Omega_t \subset \R^d$ denote a bounded open domain occupied by the solid at time $t$. Let us assume that 
there are $n+1$ different atomic species (with $n\in \N^*$) composing the material, whose local concentrations are respectively denoted by $u_0(t,x), u_1(t,x), \cdots, u_n(t,x)$ for $t\in (0,T)$ and $x\in \Omega_t$. 

The system of equations we consider here reads:
\begin{equation}\label{eq:diff}
\partial_t u_i - \mbox{\rm div}_x\left[ \sum_{0\leq j \neq i \leq n} K_{ij}(u_j \nabla_x u_i - u_i \nabla_x u_j) \right] = 0 \mbox{ in }\Omega_t, \; t\in(0,T), \; \forall \;0 \leq i \leq n,
\end{equation}
where for all $0\leq i\neq j \leq n$, the non-negative real number $K_{ij} = K_{ji} \geq 0$ denotes the cross-diffusion coefficient of atoms of type $i$ with atoms of type $j$. 
This system of equations can be formally derived as the continuous limit of some discrete stochastic lattice hopping model. We do not explain the derivation of these equations here 
for the sake of brievity and refer the reader to~\cite{Bakhta,Burger,Jungel} for more details.

The initial condition $(u_0^0, u_1^0, \cdots, u_n^0)\in L^1(\Omega_0 ; \R^{n+1})$ of this system is assumed from now on to satisfy the constraints
\begin{equation}\label{eq:initconstraint}
\forall\; 0\leq i \leq n, \; \mbox{ for almost all } x\in \Omega_0,\quad u^0_i(x)\geq 0, \quad \sum_{i=0}^n u^0_i(x) = 1 \mbox{ and } u_i(0,x) = u_i^0(x). 
\end{equation}
The last constraint $\sum_{i=0}^n u^0_i(x) = 1$ is a natural volumic constraint which encodes the fact that each site of the crystalline lattice of the solid has to be occupied 
(vacancies being treated as a particular type of atomic species).

\medskip

From a physical point of view, it is naturally expected that the solutions $(u_0, \cdots, u_n)$ of 
such a system satisfy similar constraints as the initial condition, i.e. for all $t\in (0,T)$, 
\begin{equation}\label{eq:constraint}
\forall \; 0\leq i \leq n, \; \forall t\in [0,T], \; \mbox{ for almost all }  x\in \Omega_t, \quad u_i(t,x)\geq 0 \mbox{ and } \sum_{i=0}^n u_i(t,x) = 1. 
\end{equation}
Indeed, the last equality can be infered at least formally from (\ref{eq:diff}) since by summing all these equations for $0\leq i \leq n$, we obtain that $\partial_t\left( \sum_{i=0}^n u_i\right) = 0$. 

\medskip

Despite their importance in chemistry or biology, it appears that the mathematical analysis of these systems, taking into account constraints (\ref{eq:constraint}), is quite 
recent~\cite{Jungel,Griepentrog,Mielke,Burger}. 
Let us first consider the case of a fixed domain, i.e. when $\Omega_t = \Omega_0$ for all 
$t\in [0,T]$, and of no flux boundary conditions
\begin{equation}\label{eq:BCnoflux}
\sum_{0\leq j \neq i \leq n} K_{ij}(u_j \nabla u_i - u_i \nabla u_j)\cdot \textbf{n} = 0 \mbox{ on }(0,T)\times \partial \Omega_0, \; \forall \; 0\leq i \leq n, 
\end{equation}
where $\textbf{n}$ denotes the outward normal unit vector to $\partial \Omega_0$. When all the coefficients $K_{ij}$ are identically equal to some constant $K>0$, 
it can be easily seen that (\ref{eq:diff}) boils down to a set of decoupled heat equations 
$$
\forall \; 0\leq i \leq n , \quad \partial_t u_i - \mbox{div}_x \left[ K \nabla_x u_i \right] = 0 \mbox{ in } (0,T)\times \Omega_0, 
$$
with Neumann boundary conditions $K \nabla_x u_i \cdot \textbf{n} = 0$ on $(0,T)\times \partial \Omega_0$ whose analysis becomes trivial. However, at least to our knowledge, the first proof of existence of global
weak solutions of (\ref{eq:diff}) satisfying constraints (\ref{eq:constraint}) with non-identical cross-diffusion coefficients is 
given in~\cite{Burger} for $n=2$ with coefficients $K_{ij}$ such that $K_{i0}>0$ for $i=1,2$ and $K_{12} = K_{21} = 0$. 
These results were later extended in~\cite{JungelZamponi2} to a general number of species 
$n\in \N^*$ with cross-diffusion coefficients satisfying $K_{i0} >0$ and $K_{ij} = 0$ for all $1\leq i \neq j \leq n$; the authors proved in addition the uniqueness of such weak solutions. 
In~\cite{JungelZamponi}, the case $n=2$ with arbitrary positive coefficients $K_{ij}>0$ is covered, though no uniqueness result is provided. 
The main difficulty of the mathematical analysis of such equations relies in the bounds (\ref{eq:constraint}), 
which are not obvious since no maximum principle can be proved for these systems in general. In all the articles mentioned above, the analysis framework used by the authors is the 
so-called \itshape boundedness by entropy method\normalfont, which was first used in~\cite{Burger} in the particular case mentioned above and then further developped by J\"ungel in~\cite{Jungel} 
to a more general theoretical setting. The main idea is to write the above system of equations as a formal gradient flow and derive estimates on the solutions $(u_0, \cdots, u_n)$ using the decay 
of some well-chosen entropy functional. The diffusivity matrix~\cite{MatthesZinsl} obtained for these systems is not a concave function of the densities though, so that standard gradient flow theory 
arguments (such as the minimizing movement method) cannot be applied in this context. In the first part of this note, we prove that for any $n\in \N^*$, in the case when $K_{ij} = K_{ji}>0$ for all $0\leq i \neq j \leq n$, the above system can be analyzed 
using the general theoretical framework introduced in~\cite{Jungel}. The proof heavily relies on ideas of~\cite{JungelZamponi}. Uniqueness of the solutions remains an open issue. 

\medskip

Interestingly, it seems that very few works deal with the mathematical analysis of such cross-diffusion systems in the case of moving boundary domains, 
despite their importance for instance in the modelisation of the CVD processes mentioned above. 
In the second part of this note, we present 
preliminary existence results for one-dimensional systems with moving boundaries, extending existence results obtained in the general theoretical framework introduced by J\"ungel~\cite{Jungel}.

\section{Existence on a fixed domain: formal gradient flow structure}\label{sec:fixed}

For the sake of completeness, let us recall here Theorem~2 of~\cite{Jungel}. 
\begin{theorem}[Theorem~2 of~\cite{Jungel}]\label{th:Jungel}
 Let $D$ be an open domain of $\R^n$ such that $D\subset (a,b)^n$ for some $a,b\in \R$. Let $A: u\in \overline{D} \mapsto A(u):=(A_{ij}(u))_{1\leq i,j\leq n} \in \R^{n\times n}$ be 
 a matrix-valued functional defined on $\overline{D}$ satisfying $A\in \cC^0(D; \R^{n\times n})$ and the following assumptions: 
 \begin{itemize}
  \item [(H1)] There exists a bounded from below convex function $h\in \cC^2(D, \R)$ such that its derivative $Dh:D \to \R^n$ is invertible on $\R^n$; 
  \item[(H2)] For all $1\leq i \leq n$, there exist $\alpha_i^*>0$ and $m_i>0$ such that for all $z = (z_1, \cdots, z_n)^T \in \R^n$ and $u=(u_1,\cdots, u_n)^T\in D$, 
  $$
  z^T D^2h(u)A(u) z \geq \sum_{i=1}^n \alpha_i(u_i)^2 z_i^2,
  $$
  where either $\alpha_i(u_i) = \alpha_i^*(u_i-a)^{m_i-1}$ or $\alpha_i(u_i) = \alpha_i^*(b-u_i)^{m_i-1}$;
  \item[(H3)] There exists $a^*>0$ such that for all $u\in D$ and $1\leq i,j \leq n$ such that $m_j>1$, it holds that $|A_{ij}(u)| \leq a^* |\alpha_j(u_j)|$. 
 \end{itemize}
Let $u^0 \in L^1(\Omega_0; \overline{D})$. Then there exists a weak solution $u$ with initial condition $u^0$ to 
\begin{equation}\label{eq:Jungeltrue}
\partial_t u = \mbox{div}_x(A(u)\nabla_x u) \mbox{ in } (0,T)\times \Omega_0, \mbox{ and } A(u)\nabla u \cdot \textbf{n} = 0 \mbox{ on }(0,T)\times \partial\Omega_0,
\end{equation}
such that for $(t,x)\in (0,T)\times \Omega_0$, $u(t,x)\in \overline{D}$ with $u\in L^2((0,T); H^1(\Omega_0,\R^n))$ and $\partial_t u \in L^2((0,T); H^{-1}(\Omega_0;\R^n))$.
\end{theorem}

\medskip

In this first section, we prove that the system of cross-diffusion equations introduced above, on a fixed domain $\Omega_0$, 
\begin{equation}\label{eq:diff0}
\partial_t u_i - \mbox{\rm div}_x\left[ \sum_{0\leq j \neq i \leq n} K_{ij}(u_j \nabla_x u_i - u_i \nabla_x u_j) \right] = 0 \mbox{ in }(0,T)\times \Omega_0, \; \forall \; 0 \leq i \leq n,
\end{equation}
and with no flux boundary conditions (\ref{eq:BCnoflux}), can be analyzed using the general theoretical framework introduced in~\cite{Jungel}. Actually, the following existence result is a direct consequence of Theorem~2 of~\cite{Jungel}: 

\medskip

\begin{corollary}\label{cor:cor1}
Let us assume that for all $0\leq i\neq j \leq n$, $K_{ij} = K_{ji} >0$ and let $u^0:=(u_0^0, \cdots, u_n^0)\in L^1(\Omega_0; \R^{n+1})$ satisfying condition (\ref{eq:initconstraint}). Then, there exists a weak solution 
$u = (u_0, \cdots, u_n)$ to (\ref{eq:diff0})-(\ref{eq:BCnoflux}) with initial condition $u^0$ satisfying $u\in L^2((0,T); H^1(\Omega_0; \R^{n+1}))$, $\partial_t u \in L^2((0,T);H^{-1}(\Omega_0; \R^{n+1}))$ and 
$$
\forall\; 0\leq i \leq n, \; \mbox{ for almost all }  (t,x)\in [0,T]\times \Omega_0, \quad u_i(t,x)\geq 0 \mbox{ and } \sum_{i=0}^n u_i(t,x) = 1. 
$$
\end{corollary}

\medskip

Indeed, if $(u_0, \cdots, u_n)$ is a solution to (\ref{eq:diff0}) with boundary conditions ({\ref{eq:BCnoflux}), 
it holds that $\partial_t \left( \sum_{i=0}^n u_i\right) = 0$ for almost all 
$(t,x)\in (0,T)\times \Omega_0$. Assuming that the initial condition $(u_0^0, \cdots, u_n^0)$ satisfies (\ref{eq:initconstraint}), 
the above system can then be reformulated equivalently, writing $u_0 := 1 - \sum_{i=1}^n u_i$, as a function of the concentrations $(u_1, \cdots, u_n)$ only as follows: for all $1\leq i \leq n$,
\begin{equation}\label{eq:system2}
\left\{
\begin{array}{l}
 \partial_t u_i - \mbox{\rm div}\left[\sum_{1\leq j \neq i \leq n} K_{ij}(u_j \nabla u_i - u_i \nabla u_j) + K_ {i0}( (1-\rho)\nabla u_i  + u_i\nabla \rho)\right]= 0 \mbox{ in } (0,T)\times \Omega_0, \\
 \left[\sum_{1\leq j \neq i \leq n} K_{ij}(u_j \nabla u_i - u_i \nabla u_j) + K_ {i0}( (1-\rho)\nabla u_i  + u_i\nabla \rho)\right] \cdot \textbf{n} = 0 \mbox{ on }(0,T)\times \partial \Omega_0,\\
\end{array}
\right .
\end{equation}
where $\rho:= \sum_{i=1}^n u_i = 1-u_0$. 

\medskip

In our case, $D:= \{ (u_1, \cdots, u_n)\in (\R_+^*)^n, \; \sum_{i=1}^n u_i < 1\}\subset (0,1)^n$. For any $u = (u_1, \cdots, u_n)\in \overline{D}$, 
we denote by $\rho:= \sum_{i=1}^n u_i$. The system (\ref{eq:system2}) introduced above can be rewritten under the form (\ref{eq:Jungeltrue}) using the notation
$u := (u_1, \cdots, u_n)^T$, 
with $A: u\in \overline{D} \mapsto A(u) := (A_{ij}(u))_{1\leq i,j \leq n} \in \R^{n\times n}$ defined as follows: for all $1\leq i,j \leq n$,
$$
A_{ii}(u) = \sum_{1\leq j\neq i \leq n} (K_{ij} - K_{i0}) u_j + K_{i0} \quad \mbox{ and } \quad A_{ij}(u) = -(K_{ij} - K_{i0}) u_i \mbox{ if } i \neq j.
$$

\medskip

Let us prove that $A$ satisfies the assumptions of Theorem~\ref{th:Jungel} with the following entropy functional, which is the same as in~\cite{Mielke,Jungel,JungelZamponi2}: 
\begin{equation}\label{eq:defh}
h: u\in \overline{D}\mapsto h(u):= \sum_{i=1}^n u_i\ln u_i + (1-\rho) \ln (1-\rho). 
\end{equation}

The function $h \in \cC^0(\overline{D}; \R) \cap \cC^2(D;\R)$ (thus is bounded on $\overline{D}$), is strictly convex on $D$, and its derivative $Dh: D \to \R^n$ is invertible. For all $u\in D$, 
$Dh(u) = (\ln u_i - \ln(1-\rho))_{1\leq i \leq n}$ and for all $w\in \R^n$, $Dh^{-1}(w) = \left( \frac{e^{w_i}}{1 + \sum_{j=1}^n e^{w_j}}\right)_{1\leq i \leq n}$. 
As a consequence, $h$ satisfies assumption (H1) of Theorem~\ref{th:Jungel}.  

Let us now prove that assumption (H2) of Theorem~\ref{th:Jungel} is satisfied with $m_i = \frac{1}{2}$ for all $1\leq i \leq n$. Thus, there will be no need to check assumption (H3) for the existence result to hold 
and Corollary~\ref{cor:cor1} will be a direct consequence of Theorem~\ref{th:Jungel}. To this aim, we follow the same strategy of proof as the one used in~\cite{JungelZamponi}. 
We are going to prove that if the coefficients $K_{ij}$ are assumed to be strictly positive then,  there exists $\alpha >0$ such that for all $u\in D$, 
\begin{equation}\label{eq:ineg}
H(u)A(u) \geq \alpha \Lambda(u),  \; \mbox{ where } H(u):= D^2h(u), \; \Lambda(u):= \mbox{diag}\left( \left(\frac{1}{u_i}\right)_{1\leq i \leq n}\right) \mbox{ and } \alpha:= \min_{0\leq i\neq j \leq n} K_{ij}. 
\end{equation}
This inequality implies (H2) with $m_i = \frac{1}{2}$ and $\alpha_i^* = \sqrt{\alpha}$ and $\alpha_i(u) = \sqrt{\alpha} u_i^{-1/2}$ for all $1\leq i \leq n$. 
Let $u\in D$. We have for all $1\leq i,j \leq n$,  
$$
H_{ii}(u) = \frac{1}{u_i} + \frac{1}{1-\rho} \mbox{ and } H_{ij}(u)= \frac{1}{1-\rho} \mbox{ if } i \neq j. 
$$
Introducing $P(u) := (P_{ij}(u))_{1\leq i,j \leq n}$ the matrix such that for all $1\leq i ,j \leq n$,  
$$
P_{ii}(u) = 1-u_i \mbox{ and } P_{ij}(u) = -u_i \mbox{ if } i \neq j, 
$$
it holds that $H(u) P(u) = \Lambda(u)$. Thus, we can write $H(u) A(u) - \alpha \Lambda(u) = H(u)(A(u) - \alpha P(u))$. We can easily check that 
$A(u) - \alpha P(u) = \widetilde{A}(u) + \alpha D(u)$, where $\widetilde{A}(u)$ has the same structure as $A(u)$ but with diffusion coefficients $K_{ij}-\alpha$ instead of $K_{ij}$, and 
$D(u) := (D_{ij}(u))_{1\leq i,j \leq n}$ where $D_{ij}(u) = u_i$ for all $1\leq i \leq n$.

On the one hand, $H(u)D(u) = \frac{1}{1-\rho}Z$ where $Z := (Z_{ij})_{1\leq i,j \leq n}$ with $Z_{ij} = 1$ for all $1\leq i,j \leq n$.
Since the matrix $Z$ is a positive matrix, so is $H(u)D(u)$. 

On the other hand, since $h$ is strictly convex on $D$, $H(u) \widetilde{A}(u)$ is positive if and only if $M(u):= \widetilde{A}(u) H(u)^{-1}$ is positive. Indeed, for all $z \in \R^n$, 
we have $z^T H(u)\widetilde{A}(u) z = (H(u) z)^T \left(\widetilde{A}(u) H(u)^{-1}\right) (H(u) z)$. Besides, it can be easily checked that  
$M(u) = (M_{ij}(u))_{1\leq i,j \leq n}$, where for all $1\leq i,j \leq n$,
$$
M_{ii}(u) = (K_{i0}-\alpha) (1-\rho) u_i + \sum_{1\leq j \neq i \leq n} (K_{ij}-\alpha) u_i u_j \mbox{ and } M_{ij}(u) = - (K_{ij}-\alpha) u_i u_j \mbox{ if } j\neq i,  
$$
and that this matrix is indeed a positive matrix. 
Hence we have proved inequality (\ref{eq:ineg}), assumption (H2) and Corollary~\ref{cor:cor1}.

\section{A simple one-dimensional model with moving boundary}

As announced in the introduction, the second part of this note deals with a model describing the evolution of the concentration of different chemical species in a solid thin film 
during a CVD process. In this case, the evolution of the boundary of the film has to be taken into account as well as the diffusion phenomena occuring in the bulk. We introduce and analyze 
a simplified one-dimensional model, considered as a preliminary step before tackling more complicated systems. 

In this section, $d = 1$ and at time $t=0$, the domain occupied by the solid is an interval $\Omega_0 = (0,e_0)$ for some initial thickness $e_0>0$. The external fluxes of volatile precursors during the deposition 
process are denoted by $(\phi_0, \cdots, \phi_n) \in L^\infty((0,T); \R^{n+1})$ and are assumed to satisfy $\phi_i(t)\geq 0$ for all $t\in (0,T)$ and $0\leq i \leq n$. At time $t\in(0,T)$, 
the domain occupied by the solid is equal to $\Omega_t = (0, e(t))$ where the depth $e(t)$ of the film is assumed to obey an evolution law characterized only by the values of the external fluxes: 
$$
e(t) = e_0 + \int_0^t \sum_{i=0}^n \phi_i(s)\,ds.
$$
The concentrations $(u_0, \cdots, u_n)$ satisfy the evolution law given in (\ref{eq:diff}) with the following boundary conditions: for all $0\leq i \leq n$ and $t\in (0,T)$, 
$$
\left(\sum_{0\leq j\neq i \leq n} K_{ij}(u_j \partial_x u_i - u_i \partial_x u_j)\right)(t, 0) = 0 \mbox{ and } \left(\sum_{0\leq j\neq i \leq n} K_{ij}(u_j \partial_x u_i - u_i \partial_x u_j)\right)(t, e(t)) + e'(t)u_i(t, e(t)) = \phi_i(t).
$$
These boundary conditions ensure that the total mass of each chemical species is conserved in the sense that for all $t\in (0,T)$, $\partial_t \left( \int_{0}^{e(t)} u_i(t,x)\,dx \right) = \phi_i(t)$.

\medskip

Actually, we can prove the following existence result in a slightly more general framework, directly inspired from the setting of~\cite{Jungel}.  
\begin{proposition}\label{prop:moving}
Let $D:=\{(u_1, \cdots, u_n)^T \in (\R_+^*)^n, \; \sum_{i=1}^n u_i < 1\} \subset (0,1)^n$ and $A: u\in \overline{D} \mapsto A(u):=(A_{ij}(u))_{1\leq i,j\leq n} \in \R^{n\times n}$ 
be a matrix-valued functional defined on $\overline{D}$ satisfying $A\in \cC^0(D; \R^{n\times n})$ and assumptions (H1)-(H2)-(H3) of Theorem~\ref{th:Jungel} with an entropy functional $h$ 
satisfying in addition
 \begin{itemize}
  \item[(H4)] $h\in L^\infty(\overline{D})$.
 \end{itemize}
Let $e_0>0$, $u^0\in L^1((0,e_0); \overline{D})$ and $(\phi_0, \cdots , \phi_n)\in L^\infty((0,T); \R_+^{n+1})$. We denote by $\phi:=(\phi_1, \cdots, \phi_n)^T$. Then there exists a weak solution $u$ with initial condition $u^0$ to 
\begin{equation}\label{eq:moving}
\left\{
\begin{array}{l}
e(t)= e_0 + \int_0^t \sum_{i=0}^n \phi_i(s)\,ds,\\
\partial_t u - \partial_x(A(u)\partial_x u) = 0 \mbox{ for } t\in (0,T) \mbox{ and }x\in(0, e(t)),\\
(A(u)\partial_x u)(t,e(t)) + e'(t) u(t,e(t))  = \phi(t) \mbox{ for }t\in (0,T), \\
(A(u)\partial_x u)(t,0) = 0 \mbox{ for } t\in (0,T). \\
 \end{array}
 \right .
\end{equation}
such that for $(t,x)\in (0,T)\times (0,e(t))$, $u(t,x)\in \overline{D}$ with $u\in L^2((0,T); H^1((0,e(t));\R^n))$ and $\partial_t u \in L^2((0,T); H^{-1}((0,e(t));\R^n))$.
\end{proposition}

\medskip

Note that the main difference between the assumptions of Proposition~\ref{prop:moving} and Theorem~\ref{th:Jungel} is the fact that the entropy function $h$ is asked in addition to belong to $L^\infty(\overline{D}; \R)$. 
Indeed, a crucial point in our proof is that this functional $h$ has to be bounded from above and from below on $\overline{D}$, which is indeed the case for the function $h$ defined in (\ref{eq:defh}). We refer 
the reader to~\cite{Bakhta} for a proof of this result. The idea of the proof is to define the problem onto a reference fixed domain and derive entropy estimates which enable to carry on the same analysis as in~\cite{Jungel}.

\section{Conclusion}\label{sec:conclusion}
In this note, we have presented preliminary existence results for a cross-diffusion model with a moving boundary domain, modeling the evolution of the concentrations of different chemical species in a thin film solid layer 
grown by CVD.

In the first section, we have slightly extended results of~\cite{JungelZamponi} in the sense that we have proved that the system of PDEs (\ref{eq:diff0}) with boundary conditions 
(\ref{eq:BCnoflux}) defined on a fixed domain, 
with coefficients $K_{ij} = K_{ji} >0$ for all $0\leq i\neq j \leq n$, fell into the analysis framework 
developped by J\"ungel in~\cite{Jungel}. In the second section, we have proposed and analyzed a simplified one-dimensional model in order to describe the growth of a 
thin film layer during a CVD process, as well as 
the diffusion of the different chemical species inside the bulk of the solid. The uniqueness of the solution of such a system of equations remains an open issue in both cases.  

\medskip

These results should be seen as preliminary steps before tackling more challenging systems of PDEs. In particular, the analysis of more realistic two-dimensional and three-dimensional models 
for CVD deposition taking into account surface diffusion 
phenomena is currently work in progress.

% etc, etc


%\appendix

%\section{Formal derivation of the model}\label{sec:derivation}

% The Acknowledgements are an un-numbered section
\section*{Acknowledgements}

We are very grateful to Eric Canc\`es and Tony Leli\`evre for very helpful advice and discussions.
We would like to thank Jean-Fran\c cois Guillemoles, Marie Jubeault, Torben Klinkert and Sana Laribi from IRDEP (Institut de Recherche et D\'eveloppement sur l'Energie Photovolta\"ique), 
who introduced us to the problem of modeling of CVD processes for the fabrication of thin film solar cells. The EDF company is acknowledged for funding.
We would also like to thank Martin Burger, Ansgar J\"ungel and Daniel Matthes for very helpful discussions on the theoretical part.

\end{document}